\documentclass{amsart}

\newtheorem{theorem}{Theorem}[section]
\newtheorem{lemma}[theorem]{Lemma}

\theoremstyle{definition}

\theoremstyle{remark}
\newtheorem{remark}[theorem]{Remark}

\numberwithin{equation}{section}

\begin{document}

\title[Strictly nearly K\"ahler 6-manifolds and symplectic forms]{\bf Strictly nearly K\"ahler 6-manifolds are not compatible with symplectic forms}

\author{Mehdi LEJMI}
\address{ D{\'e}partement de Math{\'e}matiques\\
UQAM\\ C.P. 8888 \\ Succ. Centre-ville \\ Montr{\'e}al (Qu{\'e}bec) \\
H3C 3P8 \\ Canada}
\email{lejmi.mehdi@courrier.uqam.ca}
\thanks{The author thanks Prof. V.~Apostolov for his  help and judicious advice, T. Draghi\u{c}i, A. Moroianu and the referee for valuable suggestions. He is very grateful to Prof. R. Bryant who
pointed out to him how some of the results in this paper are
related to \cite{bryant-1}.}

\maketitle

\begin{abstract}
We show that the almost complex structure underlying a
non-K\"ahler,  nearly K\"ahler $6$-manifold (in particular, the
standard almost complex structure of $S^6$) cannot be compatible
with any symplectic form, even locally.
\end{abstract}

\section{introduction}

Every symplectic manifold $(M, \omega)$ gives rise to an infinite
dimensional, contractible Fr\'echet space of {\it
$\omega$-compatible} almost complex structures, $J$,  introduced
by the property that the bilinear form $g(\cdot, \cdot)= \omega
(\cdot, J \cdot)$ is symmetric and positive-definite (i.e. defines
a Riemannian metric on $M$); in this case the $(J, g)$ is an
almost Hermitian structure on $M$,  which is referred to as an
{\it almost K\"ahler} structure compatible with $\omega$.

It is natural to wonder whether or not a given almost complex
structure $J$ on $M$ is $\omega$-compatible for some symplectic
form $\omega$? This question, which was first raised and studied
by J.~Armstrong in \cite{armstrong}, can be asked both locally and
globally and  the corresponding answers are quite different in
nature. In  this Note we are interested in the local aspect of the
problem, namely we consider the following

\vspace{0.1cm}\noindent {\it Question  1.} Is a given almost
complex structure $J$ on $M$ locally compatible with symplectic
forms? In other words, given $J$, can one find in a neighbourhood
of each point of $M$ a symplectic form compatible with $J$?

\vspace{0.1cm} Following~\cite{riviere-tian-2}, we shall refer to
almost complex structures which are locally compatible with
symplectic forms as almost complex structures having the {\it
local symplectic property}.  Further motivation for studying in a greater detail  this property comes from~\cite{riviere-tian-1,riviere-tian-2}, where it shown that many features of the theory of pseudo-holomorphic mappings and currents on symplectic and complex manifolds can be extended to compact almost complex manifolds having the local symplectic property.

As a trivial example, any {\it integrable} almost complex
structure satisfies the local symplectic property (although there
are many complex manifolds which are not symplectic). In
particular, when $M$ is $2$-dimensional the answer of Question~1
is always positive.

\section{The 4-dimensional case}
We start by providing a detailed proof of the following observation made in \cite[p.~10]{armstrong}:
\begin{theorem}\label{th1}
Any almost complex manifold of dimension $4$ has the local
symplectic property.
\end{theorem}
A proof of this result is readily available in
\cite[Lemma~A.1]{riviere-tian-1}. For the sake of completeness,
and since we find the arguments in \cite{riviere-tian-1}
incomplete (see Remark~\ref{error1} below), we give here an
alternative argument  based on the Malgrange existence theorem of
local solutions of elliptic systems of PDE's.

\vspace{0.1cm}\noindent {\it Proof of Theorem~\ref{th1}.} Let
$(M,J)$ an almost complex $4$-manifold.  The vector bundle of
(real) $2$-forms, $\wedge ^{2}(M)$,  decomposes with respect to
$J$ as a direct sum $\wedge^2(M) = \wedge^{J,+}(M) \oplus
\wedge^{J,-}(M),$ where $\wedge^{J,+}(M)$
(resp.~$\wedge^{J,-}(M)$) is the vector bundle of $\wedge
J$-invariant 2-forms (resp.~$\wedge J$-anti-invariant) 2-forms.
(The vector bundle $\wedge^{J,-}(M)$, endowed with the complex
structure $({\mathcal J} \phi) (\cdot, \cdot) : = -\phi(J\cdot,
\cdot),$ is naturally isomorphic to the anti-canonical bundle
$K^{-1}_J \cong \wedge^{0,2}(M)$ of $(M,J)$; likewise,
$\wedge^{J,+}(M)\otimes{\bf C} \cong \wedge^{1,1}(M)$.) We denote
by $\Omega^{J,\pm}(M)$ etc. the spaces of smooth sections  of the
corresponding bundles. The above splitting of real $2$-forms gives
rise to a decomposition of the exterior derivative $d\colon \Omega
^{1}(M)\rightarrow \Omega ^{2}(M)$ as the sum of two differential
operators $d^{\pm}  \colon \Omega ^{1}(M)\rightarrow \Omega ^{J,
\pm}(M)$.

In order to prove Theorem~\ref{th1}, it is enough to
show that  for any point $p \in M$ there exists a (connected)
neighborhood $U \ni p$ and a  $1$-form  $\alpha \in
\Omega^{1}(U)$, such that, on $U$,
\begin{equation}\label{conditions}
d^{-} \alpha=0,  \ \   d \alpha \wedge d \alpha > 0,
\end{equation}
where the sign of a $4$-form  is determined by the orientation
induced by $J$. Indeed, the $2$-form  $\omega = d \alpha = d^+\alpha$ will be then symplectic and $\wedge J$-invariant. It follows that at each point of $U$, $g(\cdot, \cdot) := \omega(\cdot, J\cdot)$ is a hermitian-symmetric 2-form on $(T(M),J)$, which can be diagonalized with respect to an Hermitian product $h$; the condition $\omega \wedge \omega >0$ means that $g$ has a positive determinant with respect to $h$; since we are in complex dimension $2$, the latter condition implies that $g$ is either positive or negative definite at the given point (and hence by continuity everywhere on $U$); the almost complex structure $J$ is, therefore,  compatible with either $\omega$ or $-\omega$.

To solve \eqref{conditions}, we first notice that the principal
symbol of $d^{-}$ is the linear map
 $$\sigma(d^-)_\xi (\alpha) = \frac{1}{2}\big(\xi \wedge \alpha -J^*\xi
\wedge J^*\alpha\big),$$ where $\xi, \alpha \in T_p^{\ast }(M)$ and
$J^*$ acts on $T_p^{\ast }(M)$ by $(J^*\alpha )(X)=-\alpha (JX)$.
Thus,  in $4$ dimensions,  $\sigma(d^-) _{\xi } : \wedge^1_p(M)
\mapsto \wedge^{J,-}_p(M)$ is {\it surjective} for any $\xi\in
T_{p}^{\ast }M\setminus \{0\}$.

 We can then associate to $d^-$ a second order {\it elliptic} linear
differential operator $P : \Omega^{J,-}(M) \to \Omega^{J,-}(M)$ by
putting $P := d^- \delta^h,$ where $h$ is some $J$-compatible
almost Hermitian metric on $(M,J)$ and $\delta^h : \Omega^2(M) \to
\Omega^1(M)$ is the corresponding co-differential, the formal
adjoint operator of $d$ with respect to $L_2$-product defined by
$h$. (The principal symbol of $P$ is given by $\sigma(P)
_{\xi}(\Phi )= - \left \vert \xi \right\vert ^{2}\Phi, \ \forall
\xi \in T^*_p(M), \Phi \in \wedge_p^{J,-}(M)$).

In terms of $P$, we want to show that for any given point $p \in
M$ one can a find a (connected) neighborhood $U \ni p$ and a
$J$-anti-invariant $2$-form $\Phi \in \Omega^{J,-}(U)$, such that
\begin{equation}\label{conditions1}
P(\Phi)=0,  \ \   d\delta^h (\Phi) \wedge d\delta^h (\Phi) > 0
\end{equation}
at any  point of $U$. Since $P$ is elliptic,  it is enough to find
 a smooth $2$-form  $\Phi_0 \in \Omega^{J,-}(M)$,  which verifies \eqref{conditions1} only at $p$ (i.e. an infinitesimal solution of \ref{conditions1}). Indeed, for any such $\Phi_0$ one can consider the system
$P(\Psi )+P(\Phi _{0})=0.$
Using the implicit function theorem, it is shown
in~\cite[p.~132]{malgrange} that for any $\varepsilon >0$ there
exist a neighborhood $U_{\varepsilon }$ of $p$ and a solution
$\Psi _{\varepsilon }\in \Omega ^{J,-}(U_{\epsilon })$  with
$\left\Vert \Psi _{\varepsilon }\right\Vert_{C^{2,\alpha }}
<\varepsilon$ (where $||\cdot ||$ stands for the H\"older norm of
$C^{2,\alpha}(U)$). Then, for $\varepsilon$ small enough, $\Phi
=\Phi _{0}+\Psi _{\varepsilon }$ and $U=U_{\varepsilon}$ will
satisfy \eqref{conditions1}.

We thus reduced the problem to verifying that at each point $p\in
M$ an infinitesimal solution always exists (for a suitable choice
of $h$). Denote by $S^{\ell}(T^*_p(M))\otimes \wedge_p^{J,-}(M)$
the space of $\ell$-jets at $p$ of elements of $\Omega^{J,-}(M)$
(where $S^{\ell}$ stands for the $\ell$-th symmetric tensor
power). By the Borel lemma, for any sequence $a_{\ell} \in
S^{\ell}(T^*_p(M))\otimes \wedge_p^{J,-}(M), \ (\ell =0,1, \cdots
)$, there exists a $\Phi \in \Omega^{J,-}(M)$ whose $\ell$-th jet
at $p$ is  $a_{\ell}$. Thus, it is enough  to show that there
exists  jets $e=(a_2,a_1,a_0)$ such that $P(e)=0$ and
$d\delta^h(e) \wedge d \delta^h (e) >0$,  where
 the linear differential operators of order $\le 2$ are identified with the induced linear maps on the space of jets of order $\le 2$. In fact, we will seek for an $e$ verifying the yet
stronger condition $((d\delta^h (e))_0=0$, where $(\cdot)_0$
denotes the primitive  part of a $2$-form  (i.e. the orthogonal
projection to $F^{\perp}$). Clearly, $(d\delta^h (e))_0=0$ implies
$P(e)=0$ and  $d\delta^h (e) =\frac{1}{2}L_h(e)F$, where $L_h$
corresponds to the the linear differential operator $L_h(\Phi)
:=h(d\delta^h \Phi, F)$. It follows that $d\delta^h (e) \wedge
d\delta^h (e) = \frac{1}{2} (L_h(e))^2 v_h$ which is positive as
soon as $L_h(e) \neq 0$. A standard calculation shows that $L_h$
is, in fact, of order one and has principle symbol $$\sigma_{\xi}
(L_h) (\Phi) = - \Phi(\xi^{\sharp}, J \theta_h^{\sharp}) + 2
\sum_{i=1}^{4} \Phi(JN(\xi^{\sharp},e_i),e_i),$$ where $\theta^h :
= J \delta_h F$ is the Lee form of $(h,J)$, $\sharp$ stands the
isomorphism between  $T^*(M)$ and $T(M)$ via $h$, $\{e_i\}$ is any
$h$-orthonormal basis of $T_p(M)$ and $$4N(\cdot, \cdot)=
\left[J\cdot,J\cdot\right] -J \left[J\cdot, \cdot\right]
-J\left[\cdot,J\cdot\right] -\left[\cdot,\cdot \right]$$ is the
Nijenhuis tensor of $J$. By making  a conformal
change $e^f h$  with $f(p)=0, df_p \neq 0$,
if necessary, we may assume that  $\sigma_{\xi} (L_h) \neq 0$. Thus, we can start
with $e'=(a_1,a_0)$ such that $L_h(e') \neq 0$. The principal
symbol of $(d \delta^h)$ is $$\sigma_{\xi}(d \delta^h)(\Phi)= - \xi
\wedge \iota_{\xi^{\sharp}}\Phi.$$ By polarization over $\xi$, it
induces a linear map from $S^2(T^*_p(M)) \otimes
\wedge_p^{J,-}(M)$ to the space of primitive $2$-forms
$(\wedge^2_p(M))_0$ which turns out  to be {\it surjective}. This
tells us that there exists an $a_2 \in S^2(T^*_p(M))\otimes
\wedge^{J,-}_p(M)$ such that $e=(a_2,a_1,a_0)$ verifies $(d
\delta^h (e))_0=0$. Since $L_h(e)=L_h(e') \neq 0$, this concludes
the proof.

\begin{remark}\label{error1} {\rm The argument given in \cite{riviere-tian-1} relies on a claim from  \cite{olver} that for {\it any} non-degenerate 2-form $\Omega$ with $d\Omega \neq 0$, there exists a local system of  coordinates $(x,y,z,t)$ such that
$$\Omega = e^x(dx\wedge dy + dz\wedge dt).$$
We note that the existence of such coordinates implies that
$\Omega$ is conformal to a symplectic form. There are, however,
many non-degenerate $2$-forms which do not verify the latter
condition. Indeed, in $4$ dimensions, to any non-degenerate 2-form
$\Omega$ one can associate a $1$-form $\theta$, called the {\it
Lee form}, such that $d\Omega = \theta \wedge \Omega$. Under a
conformal transformation ${\widetilde \Omega} = e^f \Omega$ the
Lee form changes by ${\widetilde \theta} = \theta + df$. It
follows that   $\Omega$ is
 (locally) conformally symplectic iff $d\theta =0$. For example, the $2$-form $\Omega =e^{xz}dx\wedge
dy+dz\wedge dt$ is non-degenerate and has non-closed Lee form
$\theta= e^{xz}dz$.}
\end{remark}

\begin{remark}\label{error2} {\rm Theorem~3.1 in  \cite{tomassini} affirms that there are almost complex structures on ${\bf R}^4$, which do {\it not} obey the local symplectic property.
One can see that the statement is incorrect by constructing
symplectic forms compatible with these almost complex structures.
In fact, when the function $f(x)$ in this theorem depends on $x_3$
only, then the corresponding almost complex structure is even
integrable. }
\end{remark}

\section{Strictly nearly
K\"ahler 6-manifolds}

The situation dramatically changes in dimension greater than $4$.
Indeed, it follows from~\cite{bryant-2} that the standard almost
complex structure of $S^6$ does not satisfy the local symplectic
property; A.~Tomassini~\cite{tomassini} gave other explicit
examples of $6$-dimensional almost complex manifolds which do not
satisfy the local symplectic property. In dimension greater than
$10$, J.~Armstrong~\cite{armstrong} proved that there is an open
set of (germs of) almost complex structures which doesn't satisfy
the local symplectic property. Nevertheless, a criterion of
deciding if a given almost complex structure has the local
symplectic property  is still to come.

\vspace{0.2cm}

We give below a negative answer to Question~1 for a special class
of almost complex 6-manifolds of increasing current interest, the
so-called {\it strictly nearly K\"ahler} 6-manifolds (see
e.g.~\cite{butruille,gray,moroianu-et-al,nagy,carrion,verbitsky}
and the references therein).\footnote{After the submission of a
first version of the manuscript, it was kindly pointed out to me
by R. Bryant that this result also follows from the more general
considerations in \cite[Sec.3]{bryant-1}.}
\begin{theorem}\label{th2}
The underlying almost complex structure of a non-integrable,
nearly K\"ahler 6-manifold is not compatible with any symplectic
form.
\end{theorem}

Recall that an almost Hermitian structure $(h,J)$ is {\it nearly
K\"ahler} if the covariant derivative (with respect to the
Levi-Civita connection $D^h$) of the corresponding fundamental
2-form $F(\cdot, \cdot) = h(J\cdot, \cdot)$ satisfies
$D^hF = \frac{1}{3}dF$ (nearly K\"ahler manifolds was first
studied by A. Gray \cite{gray}).
Equivalently, the {Nijenhuis thensor} $N$ is related to $dF$ by
(see e.g. \cite{kobayashi-nomizu}):
\begin{equation}\label{NK-N}
h(JN(X,Y),Z) = \frac{1}{3}dF(X,Y,Z), \ \ \forall X,Y,Z \in T(M).
\end{equation}
Apart from the integrable case, examples include $S^6$ with its
canonical almost complex structure and metric, the bi-invariant
almost complex structure on $S^3\times S^3$ with its 3-symmetric
almost-Hermitian structure, the twistor spaces over Einstein
self-dual $4$-manifolds, endowed with the anti-tautological almost
complex structure.

A key property of a {\it non-integrable} nearly K\"ahler
6-manifold is that the 3-form $dF$ is the imaginary part of a
nowhere vanishing complex $(3,0)$-form $\Psi$ on
$(M,J)$~\cite{carrion}. The identity \eqref{NK-N} then reads as
\begin{equation}\label{nijenhuis}
N =\frac{1}{6} h^* \circ \Psi,
\end{equation}
where the Nijenhuis tensor $N$ is viewed as a linear map $N :
\wedge^2(T^{1,0}(M)) \to T^{0,1}(M)$, the induced Hermitian metric
$h^*$ on $(T^*M, J^*)$ provides an isomorphism $h^* :
\wedge^{1,0}(M) \to T^{0,1}(M)$, and the complex volume form
$\Psi$ identifies $\wedge^2(T^{1,0}(M))$ with $\wedge^{1,0}(M)$.

Theorem~\ref{th2} is then an immediate corollary of the following
\begin{lemma}~\label{main}
 Let $(M,J)$ be an almost complex 6-manifold. Suppose that at some point $p$ the Nijenhuis tensor $N$ does not vanish and can be written in the form
\begin{equation}\label{algebraic}
N_p = h^*_p \circ \psi_p,
\end{equation}
where $h^*_p : \wedge^{1,0}_p(M) \to T^{0,1}_p(M)$ defines a real,
$J^*_p$-invariant, symmetric quasi-definite  form  on $T^*_p(M)$,
and $\psi_p\in \wedge^{3,0}_p(M)$ is a non-zero $(3,0)$-form.
Then, $J$ cannot be compatible with any symplectic form defined in
a neighborhood of $p$.
\end{lemma}
\noindent {\it Proof of Lemma~\ref{main}.} Since $h^*_p$ is
$J^*_p$-invariant, symmetric and quasi-definite,   there exists a
basis  $\{ \alpha_1, \alpha_2, \alpha_3 \}$ of
$\wedge^{1,0}_p(M)$,  with dual basis $\{Z_1, Z_2, Z_3\}$ of
$T^{1,0}_p(M)$, such that $h^*_p = \sum_{i=1}^3 \lambda_i (Z_i
\otimes {\overline Z_i} + {\overline Z_i} \otimes Z_i)$ with
$\lambda_i \ge 0$,  and $\psi_p = \alpha_1 \wedge \alpha_2 \wedge
\alpha_3$. (Since $N_p$ is not zero, at least one of the
$\lambda_i$'s is positive.) The condition \eqref{algebraic} then
reads as
\begin{equation}\label{Local structure}
N(Z_{1},Z_{2})=\lambda _{3}\overline{Z}_{3}, \ \ N(Z_{2},Z_{3})=\lambda _{1}%
\overline{Z}_{1}, \ \ N(Z_{3},Z_{1})=\lambda _{2}\overline{Z}_{2}.
\end{equation}

Suppose $J$ is $\omega$-compatible for some symplectic form about
$p$. The corresponding  almost K\"{a}hler structure $(J,g,\omega
)$ then satisfies (see e.g. \cite{kobayashi-nomizu}):
$(D^{g}_X\omega)(Y,Z)=-2g(JN(Y,Z),X),$ where $D^g$ is the
Levi-Civita connection of $g$. Taking a cyclic permutation  over
$X,Y,Z$ and using the fact that $\omega$ is closed, one gets
$$\underset{X,Y,Z}{\sigma }(g(JN(Y,Z),X))=0.$$ With respect to the
local basis verifying (\ref{Local structure}), the latter equality implies
$$\sqrt{-1}\sum\limits_{i=1}^{3}\lambda _{i}||Z_i||_g^2=0,$$ a
contradiction.
\begin{remark} \label{generalization} {\rm The proof of Lemma~\ref{main} shows slightly more: there is no an almost Hermitian metric $g$,  defined in a neighborhood of $p$, such that the fundamental 2-form $\omega$ of $(g,J)$ satisfies $(d\omega)^{3,0} =0$, where $(d\omega)^{3,0}$ stands for the projection of $d\omega$ to $\wedge^{3,0}(M)$.}
\end{remark}

\end{document}